# A Theoretical Study of Mafia Games

Erlin Yao[1]

Key Lab of Intelligent Information Processing,

Institute of Computing Technology, Chinese Academy of Sciences, 100190 Beijing, China

## Abstract

Mafia can be described as an experiment in human psychology and mass hysteria, or as a game between informed minority and uninformed majority. Focus on a very restricted setting, Mossel et al. [to appear in Ann. Appl. Probab. Volume 18, Number 2] showed that in the mafia game without detectives, if the civilians and mafias both adopt the optimal randomized strategy, then the two groups have comparable probabilities of winning exactly when the total player size is $R$ and the mafia size is of order $\sqrt{R}$. They also proposed a conjecture which stated that this phenomenon should be valid in a more extensive framework. In this paper, we first indicate that the main theorem given by Mossel et al. [to appear in Ann. Appl. Probab. Volume 18, Number 2] can not guarantee their conclusion, i.e., the two groups have comparable winning probabilities when the mafia size is of order $\sqrt{R}$. Then we give a theorem which validates the correctness of their conclusion. In the last, by proving the conjecture proposed by Mossel et al. [to appear in Ann. Appl. Probab. Volume 18, Number 2], we generalize the phenomenon to a more extensive framework, of which the mafia game without detectives is only a special case.

[1] Corresponding address: Institute of Computing Technology (ICT), Chinese Academy of Sciences, 100190 Beijing, China. Tel.: +86 10 58813136.
E-mail address: yaoerlin@hotmail.com.

## 1. Introduction

Mafia can be described as an experiment in human psychology and mass hysteria, or as a game between informed minority and uninformed majority. The history of mafia game can be traced to the 1970's, when a game called Murder, which was the prototype of mafia game, had been played for many years [1]. This is a game which can train one's expressive ability, discrimination ability, imagination ability and acting ability, etc. The game has common characteristics with many real life games, for example: workers, managers and stockholders at a company; students, teachers and management of a school; or citizens, the mafia and the police at a certain city. All of these games share the following three similar features: (1) each player may belong to one or more coalitions; (2) different groups make decision in different ways and take actions of different types; (3) different players accumulate information in different ways [2, 3].

According to the classification of games, the mafia game is a type of dynamic, partial information and group game. However, the previous works on the partial information and group games are mostly concerned with general definitions and abstract results in the context of extensive games. This line of research has not resulted in much quantitative analysis [3, 4]. Focus on a very restricted setting, Braverman, Etesami and Mossel obtained very precise results on the relative power of different groups. In particular, they first find a randomized strategy that is optimal in the absence of detectives, which leads to a stochastic asymptotic analysis using martingale arguments and the martingale stopping time theorem [5], where it is shown that the two groups have comparable winning probabilities when the mafia size is of order $\sqrt{R}$ . Then they propose a conjecture which states that the following phenomena should be valid in further generality: "In cases where there exists a distinguished group of size $M$ that has complete information and acts at all rounds playing against a group of players of size $R-M$ with no prior information that acts only at fraction $\alpha$ of the rounds, it is expected that the two groups will have comparable winning probabilities if $M = R^{\alpha}$" [2].

This paper investigates the quantitative analysis of mafia game without detectives and adopts the same setting of mafia game as [2]. We first note that the main theorem given by Braverman, Etesami and Mossel can not hold under the conditions given by them. More important, by constructing an example, we indicate that even the holding of the theorem can not guarantee their conclusion, i.e., the two groups have comparable winning probabilities when the mafia size is of order $\sqrt{R}$. Then through the establishment of a recursive formula, we build almost tight bounds of the winning probability and validate the conclusion of Braverman et al. In the last, by the proof of the conjecture proposed by Braverman, Etesami and Mossel, we generalize that this phenomenon is valid in a more extensive framework, of which the mafia game without detectives is only a special case.

Section 2 introduces the main results developed by Mossel et al. and indicates the deficiencies of their works. In Section 3, by the establishment of a recursive formula of the winning probability, we derive almost tight bounds of winning probability using the technique of mathematical induction, which validates the phenomenon that in the mafia game without detectives, the citizens and mafias have comparable winning probabilities when the mafia size is of order square root of the size of the total population. In Section 4, by the proof of a conjecture proposed by Mossel et al., we generalize the mafia game to a more extensive framework. Section 5 studies the change the existence of detectives can bring to the qualitative behavior of the game, where the deficiencies of existed works and future research directions are discussed. Section 6 summarizes this paper.

## 2. Existed Works and Deficiencies

In the absence of detectives, Braverman, Etesami and Mossel found that the randomized strategy is optimal for both the citizens and mafias [2]. Given the optimal strategies, suppose initially there are in all $R$ players and $M$ mafia members among the players, and let $w(R, M)$ denote the probability that the mafias win the mafia game without detectives when both the citizens and mafias adopt their optimal strategies. By using the martingale argument in stochastic asymptotic analysis [5], Braverman et al. draw one of their main results as the following theorem [2]:

**Conclusion 1.** There exists functions $p : (0, \infty) \to (0, 1)$ and $q : (0, \infty) \to (0, 1)$ such that if $0 < \eta < \infty$, the number of residents $R$ is sufficiently large and the mafia size satisfies $M \in [\eta\sqrt{R}, \eta\sqrt{R} + 1]$ then

$$p(\eta) \le w(R, M) \le q(\eta).$$

Furthermore,

$$\lim_{\eta\to\infty} p(\eta) = 1,$$

and

$$\lim_{\eta\to 0} q(\eta) = 0.$$

And according to the above theorem, the authors draw the conclusion that when there are no detectives the mafias and citizens have comparable chance to win when the mafia size $M$ is of order $\sqrt{R}$. Moreover if $M$ is a large multiple of $\sqrt{R}$ then the chance that the mafias win is close to 1 and if it is a small multiple of $\sqrt{R}$ then the chance that the mafias win is close to 0.

However, we note that the theorem does not hold under the conditions given above, the authors may have thought some conditions for granted during their proof process and the conditions should be slightly changed to make the theorem holds. More important, it will be indicated that according to the theorem, we can not derive the conclusion that the two groups have comparable winning probability when the mafia size $M$ is of order $\sqrt{R}$.

**Claim 1.** The function $q$: $(0, \infty) \to (0, 1)$ in conclusion 1 does not exist.

**Proof.** We note that for any finite $R$, the function $q$ does not exist. For example, if $\eta \ge \sqrt{R}$, then $M \ge \eta\sqrt{R} \ge R$, which means that all the players are mafias, so it is clear that $w(R, M)$=1 when $\eta \geq \sqrt{R}$. Suppose there exists $q$: (0, ∞) → (0, 1) such that $w(R, M) \leq q(\eta)$, then when $\eta \geq \sqrt{R}$,

$$q(\eta) \geq w(R, M) = 1,$$

this contradicts with $q(\eta)<1$, so the function $q$ does not exist. It is worth noting that even if we let $q$: (0, ∞) → (0, 1], for any finite $R$, there still does not exist function $q$ such that $w(R, M) \leq q(\eta)$ and $\lim_{\eta \to 0} q(\eta) = 0$. This is because when $\eta<1/\sqrt{R}$, then $\eta\sqrt{R}<1$, according to $M \in [\eta\sqrt{R}, \eta\sqrt{R}+1]$ and $M$ is integer, we have $M$ =1. So we have when $\eta<1/\sqrt{R}$,

$$w(R, 1) \leq q(\eta),$$

then

$$\lim_{\eta \to 0} q(\eta) \geq w(R,1).$$

It is clear that $w(R, 1)>0$ for any finite $R$, then we have

$$\lim_{\eta \to 0} q(\eta) > 0,$$

this contradicts with

$$\lim_{\eta \to 0} q(\eta) = 0. \square$$

It is clear that

$$\lim_{R \to \infty} w(R,1) = 0,$$

so $\lim_{\eta \to 0} q(\eta) = 0$ should be changed to:

$$\lim_{R \to \infty} \lim_{\eta \to 0} q(\eta) = 0.$$

Then conclusion 1 should be slightly changed as following:

**Conclusion 1'.** There exists functions $p$ : (0, ∞) → (0, 1) and $q$ : (0, ∞) → (0, 1] such that if $0 < \eta < \infty$, the number of residents $R$ is sufficiently large and the mafia size satisfies $M \in [\eta\sqrt{R}, \eta\sqrt{R}+1]$ then

$$p(\eta) \leq w(R, M) \leq q(\eta).$$

Furthermore,

$$\lim_{\eta\to\infty} p(\eta) = 1,$$

and

$$\lim_{R\to\infty}\lim_{\eta\to 0} q(\eta) = 0.$$

However, we note that according to the holding of the above conclusion about winning probability, we can not arrive at the conclusion that the two groups have comparable winning probability when the mafia size $M$ is of order $\sqrt{R}$. We prove this point by constructing the following example:

**Example 1.** Suppose

$$w(R, M)=M/(R+M),$$

we construct functions $p$: $(0, \infty) \to (0, 1)$ and $q$: $(0, \infty) \to (0, 1]$ as follows: let

$$p(\eta)=\eta/(\eta+\sqrt{R}).$$

If $\eta<\sqrt{R}/2-1/\sqrt{R}$, let

$$q(\eta)=(\eta+1/\sqrt{R})/(\eta+\sqrt{R});$$

and if $\eta\ge\sqrt{R}/2-1/\sqrt{R}$, let

$$q(\eta)=1.$$

Because $R>0$, it is clear that for any positive $x$ and $y$, if $x\ge y$, then

$$x/(R+x)\ge y/(R+y).$$

Since $\eta\sqrt{R} \le M \le \eta\sqrt{R} + 1$, then we have:

$$\eta\sqrt{R}/(R+\eta\sqrt{R}) \le M/(R+M) \le (\eta\sqrt{R} + 1)/(R+\eta\sqrt{R} + 1) \qquad (1)$$

In Eq. (1),

$$\eta\sqrt{R}/(R+\eta\sqrt{R})=\eta/(\eta+\sqrt{R})= p(\eta);$$

$$(\eta\sqrt{R} + 1)/(R+\eta\sqrt{R} + 1)= (\eta+1/\sqrt{R})/(\eta+\sqrt{R}+1/\sqrt{R})$$

$$< (\eta+1/\sqrt{R})/(\eta+\sqrt{R})= q(\eta).$$

Then according to Eq. (1) and $M/(R+M) = w(R, M)$, we have

$$p(\eta) \le w(R, M) \le q(\eta).$$

And it is clear that

$$\lim_{\eta\to\infty} p(\eta) = \lim_{\eta\to\infty} \frac{\eta}{\eta+\sqrt{R}} = 1$$

and

$$\lim_{R\to\infty}\lim_{\eta\to 0} q(\eta) = \lim_{R\to\infty}\lim_{\eta\to 0} \frac{\eta + 1/\sqrt{R}}{\eta+\sqrt{R}} = \lim_{R\to\infty} \frac{1}{R} = 0 . \square$$

It can be seen that $w(R, M)=M/(R+M)$ and the functions $p$ and $q$ satisfy all the conditions in conclusion 1'. But it is clear that in this case only when $M$ and $R$ are in the same order, $w(R, M)$ can be a comparable value with 1/2. So the holding of conclusion 1' can not guarantee that the two groups have comparable winning probability when the mafia size $M$ is of order $\sqrt{R}$ . Then what conditions should $R$ and $M$ satisfy to keep the two groups have comparable winning probabilities? The following section will answer this problem.

## 3. Almost Tight Bounds of Winning Probability

Suppose at the beginning of the game, there are in all $n$ civilians and $m$ mafias (here $m$ and $n+m$ correspond to the $M$ and $R$ in the above section respectively), let $W(n, m)$ be the probability that the mafias win the game in the last. We note that in the mafia game without detectives, the randomized strategy is optimal for both sides only when the civilians are in majority of the players. When the mafias are in majority of the players, then no matter what strategy the civilians adopt, the mafias will surely win the game, which means that $W(n, m)=1$ if $n<m$. And it is clear that $W(n, 0)=0$ for any $n>0$. The left boundary case is $n=1$ and $m=1$, it is clear that in this case there will be a tie, since the civilian and mafia will vote each other. According to the rule adopted by Braverman, Etesami and Mossel, "In cases of a tie, the identity of the person to be killed is chosen uniformly at random among all players who received the maximal number of votes" [2], we should let $W(1, 1)=0.5$.

If the civilians are in majority and both sides adopt their optimal strategy, then we have a recursive formula of $W(n, m)$ as the following theorem:

**Theorem 1.** In the mafia game without detectives, suppose initially there are $n$

civilians and $m$ mafias, the civilians are in majority and both sides adopt their optimal randomized strategy. Let $W(n, m)$ be the probability that the mafias win the game in the last, then we have:

$$W(n,m) = \frac{n}{n+m} W(n-2,m) + \frac{m}{n+m} W(n-1,m-1) \quad (2)$$

holds for any $n \geq m \geq 1$.

**Proof.** The game consists of the iteration of the day round and the night round until the game terminates. If both sides adopt the randomized strategy, then in the day round every player has the equal chance to be voted to death. So in the vote round, one of the civilians is selected with probability $n/(n+m)$ and one of the mafias is selected with probability $m/(n+m)$. In the first case, in the night round one civilian is randomly selected to kill, so after the first iteration, there are $n-2$ civilians and $m$ mafias left. In the second case, in the night round one civilian is randomly selected to kill, so after the first iteration, there are $n-1$ civilians and $m-1$ mafias left. So we have the recursive formula of $W(n, m)$ as Eq. (2). □

It would be nice if we can derive the closed-form expression of $W(n, m)$ according to Eq. (2). However, we note that the closed-form expression of $W(n, m)$ is difficult to derive even if when $m$ is very small. But we can get almost tight bounds of $W(n, m)$ as the following theorem:

**Theorem 2.** In the mafia game without detectives, suppose initially there are $n$ civilians and $m$ mafias, the civilians are in majority and both sides adopt their optimal randomized strategy. Let $W(n, m)$ be the probability that the mafias win the game in the last, if $m \leq k$ and $n+m \leq R$, then we have for any $n \geq m$:

$$\frac{\sqrt{2k-2}}{k} \frac{m}{\sqrt{n+m}} \leq W(n,m) \leq R^{\frac{1}{100}} \frac{m}{\sqrt{n+m}}. \quad (3)$$

**Proof.** Since $W(m-2, m)=W(m-1, m)=1$, it is clear that $W(n, m) \geq \frac{\sqrt{2k-2}}{k} \frac{m}{\sqrt{n+m}}$ holds when $n=m-2$ and $n=m-1$ for any $1 \leq m \leq k$. According to $W(n, 0)=0$, we have $W(n, m) \geq \frac{\sqrt{2k-2}}{k} \frac{m}{\sqrt{n+m}}$ holds when $m=0$ for any $1 \leq n$. Suppose for some $m \geq 0$, $W(n, m) \geq \frac{\sqrt{2k-2}}{k}\frac{m}{\sqrt{n+m}}$ holds for any $m-2 \leq n$. Then in the case $m+1$ ($m+1 \leq k$), we already have $W(n, m+1) \geq \frac{\sqrt{2k-2}}{k}\frac{m}{\sqrt{n+m}}$ holds for $n = m-1$ and $n = m$, suppose for some $t \geq m$, $W(n, m+1) \geq \frac{\sqrt{2k-2}}{k}\frac{m}{\sqrt{n+m}}$ holds for any $n \leq t$, then we have:

$$\begin{aligned} W(t+1, m+1) &= \frac{t+1}{t+m+2}W(t-1, m+1) + \frac{m+1}{t+m+2}W(t, m) \\ &\geq \frac{t+1}{t+m+2} \times \frac{\sqrt{2k-2}}{k}\frac{m+1}{\sqrt{t+m}} + \frac{m+1}{t+m+2} \times \frac{\sqrt{2k-2}}{k}\frac{m}{\sqrt{t+m}} \\ &= \frac{\sqrt{2k-2}}{k} \times \frac{(m+1)(t+m+1)}{(t+m+2)\sqrt{t+m}} \end{aligned} \tag{4}$$

It is obvious that

$$(1 - \frac{1}{t+m+2})^2 \geq 1 - \frac{2}{t+m+2} \tag{5}$$

Eq. (5) is equivalent to $(\frac{t+m+1}{t+m+2})^2 \geq \frac{t+m}{t+m+2}$, square root both sides, we have $\frac{t+m+1}{t+m+2} \geq \frac{\sqrt{t+m}}{\sqrt{t+m+2}}$, divide both sides by $\sqrt{t+m}$ and multiply both sides by $\frac{\sqrt{2k-2}}{k}(m+1)$, we have:

$$\frac{\sqrt{2k-2}}{k} \times \frac{(m+1)(t+m+1)}{(t+m+2)\sqrt{t+m}} \geq \frac{\sqrt{2k-2}}{k} \times \frac{(m+1)}{\sqrt{t+m+2}} \tag{6}$$

According to Eq. (4) and (6), $W(t+1, m+1) \geq \frac{\sqrt{2k-2}}{k} \times \frac{(m+1)}{\sqrt{t+m+2}}$ holds. According to the mathematical induction, we have $W(n, m+1) \geq \frac{\sqrt{2k-2}}{k}\frac{m+1}{\sqrt{n+m+1}}$ holds for any $n \geq m-1$. Since from assumption, we have $W(n, m) \geq \frac{\sqrt{2k-2}}{k}\frac{m}{\sqrt{n+m}}$ holds for any $m-2 \leq n$, again according to the mathematical induction, we can draw that

$$W(n, m) \geq \frac{\sqrt{2k-2}}{k}\frac{m}{\sqrt{n+m}} \tag{7}$$

holds for any $m \leq k$ and $n \geq m-2$.

Since $W(m-2, m)=W(m-1, m)=1$ and $(n+m)^{\frac{1}{100}} \geq 1$, it is clear that $W(n, m) \leq (n+m)^{\frac{1}{100}}\frac{m}{\sqrt{n+m}}$ holds when $n=m-2$ and $n=m-1$ for any $1 \leq m$. According to $W(n, 0)=0$, we have $W(n, m) \leq (n+m)^{\frac{1}{100}}\frac{m}{\sqrt{n+m}}$ holds when $m=0$ for any $1 \leq n$. Through the computing using electronic computer, we can show that $W(n, m) \leq (n+m)^{\frac{1}{100}}\frac{m}{\sqrt{n+m}}$ holds for any $n+m \leq 100$. Suppose for some $m \geq 0$, $W(n, m) \leq (n+m)^{\frac{1}{100}}\frac{m}{\sqrt{n+m}}$ holds for any $m-2 \leq n$. Then in the case $m+1$, we already have $W(n, m+1) \leq (n+m+1)^{\frac{1}{100}}\frac{m+1}{\sqrt{n+m+1}}$ holds for $n=m-1$, $m$ and $n+m+1 \leq 100$, suppose for some $t \geq \max\{m, 99-m\}$, $W(n, m+1) \leq (n+m+1)^{\frac{1}{100}}\frac{m+1}{\sqrt{n+m+1}}$ holds for any $n \leq t$, then we have:

$$\begin{aligned} W(t+1, m+1) &= \frac{t+1}{t+m+2}W(t-1, m+1) + \frac{m+1}{t+m+2}W(t, m) \\ &\leq \frac{t+1}{t+m+2} \times (t+m)^{\frac{1}{100}}\frac{m+1}{\sqrt{t+m}} + \frac{m+1}{t+m+2} \times (t+m)^{\frac{1}{100}}\frac{m}{\sqrt{t+m}} \\ &= (t+m)^{\frac{1}{100}} \times \frac{(m+1)(t+m+1)}{(t+m+2)\sqrt{t+m}} \end{aligned} \tag{8}$$

To continue the proof, we first introduce the following inequality:

$$1-\frac{1}{n} \leq (1-\frac{2}{n})^{\frac{1}{2}-\frac{1}{n}} \tag{9}$$

holds for any integer $n>2$. Let $f(x)=(1-2x)^{\frac{1}{2}-x}+x-1$ ($0 \leq x < 0.5$), it is obvious that $f(0)=0$. And $f'(x)=1-(1-2x)^{\frac{1}{2}-x}(1+\ln(1-2x))$, it is clear that if $0 < x < 0.5$, then $f'(x)>0$, so $f(x)> f(0)=0$ ($0< x <0.5$), which is equivalent to $1-x<(1-2x)^{\frac{1}{2}-x}$.

For integer $n>2$, let $x=1/n$, then we have $1-\frac{1}{n}\le(1-\frac{2}{n})^{\frac{1}{2}-\frac{1}{n}}$ holds.

Since $t+m+2>100$, in Eq. (9), let $n= t+m+2$, we have:

$$1-\frac{1}{t+m+2}\le(1-\frac{2}{t+m+2})^{\frac{1}{2}-\frac{1}{t+m+2}} \tag{10}$$

Since $\frac{1}{2}-\frac{1}{t+m+2}\ge\frac{1}{2}-\frac{1}{100}$ and $1-\frac{2}{t+m+2}<1$, then we have $(1-\frac{2}{t+m+2})^{\frac{1}{2}-\frac{1}{t+m+2}}\le(1-\frac{2}{t+m+2})^{\frac{1}{2}-\frac{1}{100}}$ holds. According to Eq. (10), we arrive at $1-\frac{1}{t+m+2}\le(1-\frac{2}{t+m+2})^{\frac{1}{2}-\frac{1}{100}}$, which under transformation is equivalent to $(t+m)^{\frac{1}{100}}\frac{t+m+1}{(t+m+2)\sqrt{t+m}}\le\frac{(t+m+2)^{\frac{1}{100}}}{\sqrt{t+m+2}}$, both sides multiply with ($m$+1), we have $(t+m)^{\frac{1}{100}}\frac{(m+1)(t+m+1)}{(t+m+2)\sqrt{t+m}}\le(t+m+2)^{\frac{1}{100}}\frac{(m+1)}{\sqrt{t+m+2}}$. According to Eq. (8), we have $W(t+1,m+1)\le(t+m+2)^{\frac{1}{100}}\frac{(m+1)}{\sqrt{t+m+2}}$ holds. According to the mathematical induction, we have $W(n, m+1)\le(n+m+1)^{\frac{1}{100}}\frac{m+1}{\sqrt{n+m+1}}$ holds for any $n\ge m-1$. Since from assumption we have $W(n, m)\le(n+m)^{\frac{1}{100}}\frac{m}{\sqrt{n+m}}$ holds for any $m-2\le n$, again according to the mathematical induction, we can draw that

$$W(n, m)\le(n+m)^{\frac{1}{100}}\frac{m}{\sqrt{n+m}}$$

holds for any $n\ge m-2$. Since $n+m\le R$, we have:

$$W(n, m)\le(n+m)^{\frac{1}{100}}\frac{m}{\sqrt{n+m}}\le R^{\frac{1}{100}}\frac{m}{\sqrt{n+m}} \tag{11}$$

Combining Eq. (7) and (11), we have if $m\le k$ and $n+m\le R$, then for any $n\ge m$,

$$\frac{\sqrt{2k-2}}{k}\frac{m}{\sqrt{n+m}} \le W(n,m) \le R^{\frac{1}{100}}\frac{m}{\sqrt{n+m}} \text{ holds.} \square$$

If the number of players $R$ is not very big, for example $R \le 10^{10}$, then $R^{\frac{1}{100}}$ is very close to 1. Note that during the proof process, it can be seen that the constant index 0.01 of $R$ can be enhanced to be very small, say, 0.0001. Then according to theorem 2, it is very clear that the two groups have comparable winning probabilities only when the mafia size $M$ is of order $\sqrt{R}$.

It is worth noting that the almost tight bounds we got here can explain the phenomena of experiments in [2]. Braverman et al. calculated the winning probability of a mafia of size $M=\eta\sqrt{R}$ as a function of $\eta$, suppose the function is $p(\eta)$, they draw the conclusion that: when $\eta$ is small, then $p(\eta)$ is almost linear; and when $\eta$ is big, $p(\eta)$ tends to 1. Theorem 2 can explicitly explain this phenomenon. According to $\frac{\sqrt{2k-2}}{k}\frac{M}{\sqrt{R}} \le W(R,M) \le R^{\frac{1}{100}}\frac{M}{\sqrt{R}}$ (in Eq. (3), let $m=M$ and $n+m=R$), if $\eta$ is small, then $\frac{\sqrt{2k-2}}{k}$ is close to 1, and the winning probability $p(\eta)$ lies between $\frac{\sqrt{2k-2}}{k}\eta$ and $\eta$, so it is almost linear; when $\eta$ is big (for example, approaches $\frac{k}{\sqrt{2k-2}}$), then $\frac{\sqrt{2k-2}}{k}\frac{M}{\sqrt{R}}$ approaches 1, since $\frac{\sqrt{2k-2}}{k}\frac{M}{\sqrt{R}} \le W(R,M)$, then it is clear that $p(\eta)$ tends to 1.

## 4. Generalization by Proof of a Conjecture

Braverman, Etesami and Mossel claimed that the following phenomenon should be valid in further generality: "In cases where there exists a distinguished group of size $M$ that has complete information and acts at all rounds playing against a group of players of size $R-M$ with no prior information that acts only at fraction $\alpha$ of the rounds, it is expected that the two groups will have comparable winning probabilities if $M = R^{\alpha}$". In particular, they propose the following conjecture [2]:

**Conjecture 1.** Consider a variant of the mafia game without detectives, where each $r$ rounds are partitioned into $d$ day rounds and $r$–$d$ night rounds. Then the two groups have comparable winning probabilities if $M=R^{d/r}$.

Note that the mafia game without detectives is only a special case where $d$=1 and $r$=2. We give the following theorem to prove this conjecture:

**Theorem 3.** Suppose there is a winning probability $W(n, m)$ ($0\le W(n, m) \le 1$), $W(n, 0)=0$ ($n>0$) and if $n< m$, $W(n, m)=1$, and satisfies the following recursive formula ($r$ and $d$ are integers and $r > d$):

$$W(n,m) = \frac{n}{n+m} W(n-r,m) + \frac{m}{n+m} W(n-r+d, m-d)\,. \qquad (12)$$

Then if $m \le C_1$ and $n+m \le C_2$, we have:

$$f(C_1)\frac{m}{(n+m)^{\frac{d}{r}}} \le W(n,m) \le g(C_2)\frac{m}{(n+m)^{\frac{d}{r}}} \qquad (13)$$

Where $f$ and $g$ are functions with finite values.

**Proof.** As the proof is very similar to that of Theorem 2, we shall only sketch it. Clearly we need only prove the following two inequalities.

The first inequality is (with contrast to Eq. (5) in Theorem 2):

$$(1-rx)^{\frac{d}{r}} \le 1-dx \quad \text{for } 0\le x < 1/r. \qquad (14)$$

Let $p(x) = (1-rx)^{\frac{d}{r}} + dx - 1$, it is clear that $p(0)=0$. And $p'(x) = d - d(1-rx)^{\frac{d}{r}-1}$, since $x < 1/r$ and $d<r$, we have $(1-rx)^{\frac{d}{r}-1} > 1$, then $p'(x)<0$, so $p(x) < p(0)=0$ ($0< x < 1/r$), which is equivalent to $(1-rx)^{\frac{d}{r}} < 1-dx$.

The second inequality is (with contrast to Eq. (9) in Theorem 2):

$$1-dx \le (1-rx)^{\frac{d}{r}-dx} \quad \text{for } 0\le x < 1/r. \qquad (15)$$

Let $q(x) = (1-rx)^{\frac{d}{r}-dx} + dx - 1$, it is clear that $q(0)=0$. And

$$q'(x) = d - d(1-rx)^{\frac{d}{r}-dx}(\ln(1-rx)+1)\,,$$

since $x < 1/r$ and $d<r$, we have $(1-rx)^{\frac{d}{r}-dx} < 1$ and $(\ln(1-rx)+1)<1$, then $q'(x)>0$, so

$q(x) > q(0)=0$ ($0< x <1/r$), which is equivalent to $(1-rx)^{\frac{d}{r}-dx} > 1-dx$. □

It is worth noting that the conditions, $W(n, 0)=0$ ($n>0$) and $W(n, m)=1$ if $n<m$, in Theorem 3 are prerequisite to the proof, because these conditions are the starting points of the mathematical induction, without the holding of the base, the conclusion may not be true. And if we do not consider the small cases (for example, $n+m \le 100$), then the coefficients of the bounds ($f(C_1)$ and $g(C_2)$) can be refined to be more accurate.

## 5. Further Discussions and Open Problems

Mossel et al. showed that even a single detective can change the qualitative behavior of the game dramatically. More formally they prove the following theorem [2]:

**Conclusion 2.** (1) Consider the game with 1 detective and mafia of size $M = \eta R < R/49$. Then for $R$ sufficiently large the probability that the mafia wins denoted $w(R, M, 1)$ satisfies $p(\eta, 1) \le w(R, M, 1) \le q(\eta, 1)$, where $0 < p(\eta, 1) < q(\eta, 1) < 1$ for all $\eta < 1/49$ and $q(\eta, 1) \to 0$ as $\eta \to 0$.

(2) Let $d \ge 1$ and consider the game with $d$ detectives and mafia of size $M = \eta R$, where $\eta<1/2$. Then for $R$ sufficiently large, the probability that the mafia wins, denoted $w(R, M, d)$ satisfies $w(R, M, d) \le q(\eta, d)$, where for each $\eta<1/2$ it holds that $\lim_{d\to\infty} q(\eta, d)= 0$.

They state that the theorem shows that even a single detective dramatically changes the citizen's team power: While in the game with no detectives a mafia of size $R^{1/2+\varepsilon}$ will surely win, as soon as there is one detective, the mafia will lose unless it is of size $\Omega(R)$.

However, we first note that the condition $M = \eta R$ in Conclusion 2 should be $M \in [\eta R, \eta R+1]$ as in Conclusion 1, else $M$ may not be an integer. Second, we show that $\lim_{\eta\to 0} q(\eta,1) = 0$ should be changed to $\lim_{R\to\infty}\lim_{\eta\to 0} q(\eta,1) = 0$.

**Claim 2.** The function $q$ in conclusion 2 does not exist.

**Proof.** We will prove that for any finite $R$, there does not exist function $q$ such that $w(R, M, 1) \le q(\eta, 1)$ for all $\eta < 1/49$ and $\lim_{\eta\to 0} q(\eta,1) = 0$. This is because when $\eta<1/R$ ($R>49$ and is finite), then $\eta R<1$, according to $M\in[\eta R, \eta R+1]$ and $M$ is integer, we have $M=1$. So we have when $\eta<1/R$,

$$w(R, 1, 1) \le q(\eta, 1),$$

then

$$\lim_{\eta\to 0} q(\eta,1) \ge w(R,1,1).$$

It is clear that $w(R, 1, 1)>0$ for any finite $R$, then we have

$$\lim_{\eta\to 0} q(\eta,1) > 0,$$

this contradicts with

$$\lim_{\eta\to 0} q(\eta,1) = 0. \square$$

It is clear that

$$\lim_{R\to\infty} w(R,1,1) = 0,$$

so $\lim_{\eta\to 0} q(\eta,1) = 0$ should be changed to:

$$\lim_{R\to\infty}\lim_{\eta\to 0} q(\eta,1) = 0.$$

Then conclusion 2 should be slightly changed as following:

**Conclusion 2'.** (1) Consider the game with 1 detective and mafia of size $M\in[\eta R, \eta R+1]$, where $\eta< R/49$. Then for $R$ sufficiently large the probability that the mafia wins denoted $w(R, M, 1)$ satisfies $p(\eta, 1) \le w(R, M, 1) \le q(\eta, 1)$, where $0 < p(\eta, 1) < q(\eta, 1) < 1$ for all $\eta < 1/49$ and $\lim_{R\to\infty}\lim_{\eta\to 0} q(\eta,1) = 0$.

(2) Let $d \ge 1$ and consider the game with $d$ detectives and mafia of size $M\in[\eta R, \eta R+1]$, where $\eta<1/2$. Then for $R$ sufficiently large, the probability that the mafia wins, denoted $w(R, M, d)$ satisfies $w(R, M, d) \le q(\eta, d)$, where for each $\eta<1/2$ it holds that $\lim_{d\to\infty} q(\eta, d)= 0$.

However we note that the holding of the first part of conclusion 2 can not guarantee that even a single detective dramatically changes the citizen's team power: as soon as there is one detective, the mafia will lose unless it is of size $\Omega(R)$. We construct the following example to prove this:

**Example 2.** Suppose

$$w(R, M, 1)=M/(\sqrt{R}+M),$$

we construct functions $p(\eta, 1)$ and $q(\eta, 1)$ as follows: let

$$p(\eta, 1)= \eta R/(\eta R+\sqrt{R}),$$

$$q(\eta, 1)=(\eta R+1)/(\eta R+1+\sqrt{R}).$$

Because $\sqrt{R}>0$, it is clear that for any positive $x$ and $y$, if $x\geq y$, then

$$x/(\sqrt{R}+x)\geq y/(\sqrt{R}+y).$$

Since $\eta R \leq M \leq \eta R + 1$, then we have:

$$\eta R/(\eta R+\sqrt{R}) \leq M/(\sqrt{R}+M) \leq (\eta R+1)/(\eta R+1+\sqrt{R}),$$

According to the definition, we have

$$p(\eta, 1) \leq w(R, M, 1) \leq q(\eta, 1).$$

And it is clear that

$$\lim_{R\to\infty}\lim_{\eta\to 0} q(\eta,1) = \lim_{R\to\infty}\lim_{\eta\to 0}\frac{\eta R+1}{\eta R+1+\sqrt{R}} = \lim_{R\to\infty}\frac{1}{1+\sqrt{R}} = 0. \square$$

It can be seen that $w(R, M, 1)= M/(\sqrt{R}+M)$ and the functions $p$ and $q$ satisfy all the conditions in part 1 of conclusion 2'. But it is clear that in this case only when $M$ and $\sqrt{R}$ are in the same order, $w(R, M, 1)$ can be a comparable value with 1/2. So the holding of part 1 of conclusion 2' can not guarantee that as soon as there is one detective, the mafia will lose unless it is of size $\Omega(R)$.

We then give the following example to show that the condition $\lim_{d\to\infty} q(\eta, d)= 0$ in the part 2 of Conclusion 2' is obvious and can not explain any phenomenon.

**Example 3.** Suppose

$$w(R, M, d)=M/(M+f(R, d)),$$

where $f(R, d)>0$ and satisfies $\lim_{d\to\infty} f(R, d)= \infty$. Let

$$q(\eta, d)=(\eta R+1)/(\eta R+1+ f(R, d)).$$

Because $f(R, d)>0$, it is clear that for any positive $x$ and $y$, if $x\geq y$, then

$$x/(f(R,d)+x) \geq y/(f(R,d)+y).$$

Since $\eta R \leq M \leq \eta R + 1$, then we have:

$$M/(M+f(R,d)) \leq (\eta R+1)/(\eta R+1+f(R,d)),$$

According to the definition, we have

$$w(R,M,d) \leq q(\eta,d).$$

And it is clear that

$$\lim_{d\to\infty} q(\eta,d) = \lim_{d\to\infty} \frac{\eta R+1}{\eta R+1+f(R,d)} = 0. \square$$

It can be seen that $w(R, M, d)=M/(M+f(R, d))$ and the function $q$ satisfy all the conditions in part 2 of conclusion 2'. But it is clear that in this case when $M$ and any $f(R, d)$ (i.e., $f(R, d)= d+\sqrt{R}$ or $f(R, d)=d^2$) are in the same order, $w(R, M, d)$ can be a comparable value with 1/2. So the holding of part 2 of conclusion 2' can not guarantee that the existence of detectives can dramatically change the citizen's team power. It would be interesting to investigate what changes the existence of detectives can take to the qualitative behavior of the game.

## 6. Conclusion

In the mafia game without detectives, when both the citizens and mafias adopt the optimal randomized strategy, we focus on the winning probabilities of different groups. By developing almost tight bounds of the winning probability, we validate the phenomenon that the two groups have comparable winning chances when the total player size is $R$ and the mafia size is of order $\sqrt{R}$. By the proof of a conjecture, we generalize the mafia game to the following extensive framework: In cases where there exists a distinguished group of size $M$ that has complete information and acts at all rounds playing against a group of players of size $R-M$ with no prior information that acts only at fraction α of the rounds, then the two groups have comparable winning probabilities if $M = R^{\alpha}$.

## Acknowledgement

This work was supported in part by the National Natural Science Foundation of China, and the National Basic Research Program of China Grant 2003CB317000, 2003CB317001.